# Improved Estimators in Simple Random Sampling When Study Variable is An Attribute

Rajesh Singh and Prayas Sharma

Department of Statistics, Banaras Hindu University, Varanasi (U.P.), India

**Abstract**

This article addresses the problem of estimating the population mean in the presence of auxiliary information when study variable itself is qualitative in nature. Bias and mean squared error (MSE) expressions of the class of estimators are derived up to the first order of approximation. The suggested estimators have been compared with the traditional estimator and several other estimators considered by Singh et al. (2010). In addition, we support this theoretical result by an empirical study to show the superiority of the constructed estimator over others.

## 1. Introduction

In the theory of sample surveys, it is usual to make use of the auxiliary information at the estimation stage in order to improve the precision or accuracy of an estimator of unknown population parameter of interest. Ratio, product and regression methods of estimation are good examples in this context. Many authors including Upadhyaya and Singh (1999), Abu-Dayyeh et al. (2003), Kadilar and Cingi (2005), Khoshnevisan et al. (2007), Singh et al. (2007), Singh et al. (2008) and Singh and Kumar (2011) suggested estimators using known population parameters of an auxiliary variable. But there may be many practical situations when auxiliary information is not available directly but is qualitative in nature, that is, auxiliary information is available in the form of an attribute. For example:

(a) The height of a person may depend on the fact that whether the person is male or female.
(b) The efficiency of a Dog may depend on the particular breed of that Dog.
(c) The yield of wheat crop produced may depend on a particular variety of wheat, etc.

In these situations by taking the advantage of point bi-serial correlation between the study variable y and the auxiliary attribute $\phi$ along with the prior knowledge of the population parameter of auxiliary attribute, the estimators of population parameter of interest can be constructed.

In many situations study variable is generally ignored not only by ratio scale variables that are essentially qualitative, or nominal scale, in nature, such as sex, race, colour, religion, nationality, geographical region, political upheavals (see Gujarati and Sangeetha (2007)). Taking into consideration the point bi-serial correlation between auxiliary attribute and study variable, several authors including Naik and Gupta (1996), Jhajj et al. (2006), Singh et al. (2007), Shabbir and Gupta (2007), Singh et al.(2008), Singh et al. (2010), Abd-Elfattah et al. (2010), and Singh and Solanki (2012) proposed improved estimators of population mean. All the others have implicitly assumed that the study variable Y is quantitative whereas the auxiliary variable is qualitative.

There may be situations when study variable itself is qualitative in nature. For example, consider U.S. presidential elections. Assume that there are two political parties, Democratic and Republican. The dependent variable here is the vote choice between two political parties. Suppose we let Y=1, if the vote is for a Democratic candidate and Y=0, if the vote is for republican candidate. Some of the variables used in the vote choice are growth rate of GDP, unemployment and inflation rates, whether the candidate is running for re-election, etc. For the present purposes, the important thing is to note that the study variable is a qualitative variable. One can think several other examples where the study variable is qualitative in nature. Thus, a family either owns a house or it does not, it has disability insurance or it does not, both husband and wife are in the labour force or only one suppose is, etc. (see Gujarati and Sangeetha (2007)). In this paper we propose estimators in which study variable itself is qualitative in nature.

Consider a sample of size n drawn by simple random sampling without replacement (SRSWOR) from a population of size N. Let $\phi_i$ and $x_i$ denote the observations on variable $\phi$ and x respectively for $i^{th}$ unit (i=1,2,3…N). $\phi_i = 1$, if $i^{th}$ unit of population possesses attribute $\phi$ and $\phi_i = 0$, otherwise. Let $A = \sum_{i=1}^{N}\phi_i$ and $a = \sum_{i=1}^{n}\phi_i$, denotes the total number of units in the population and sample possessing attribute $\phi$ respectively, $P = \frac{A}{N}$ and $p = \frac{a}{n}$, denotes the proportion of units in the population and sample, respectively, possessing attribute $\phi$.

Let us define,

$$e_p = \frac{(p-P)}{P}, \qquad e_1 = \frac{(\bar{x}-\bar{X})}{\bar{X}}, \qquad e_3 = \frac{(s_x^2 - S_x^2)}{S_x^2}$$

Such that,

$$E(e_i) = 0, (i = p, 1, 3)$$

and

$$E(e_p^2) = fC_p^2, \qquad E(e_1^2) = fC_x^2, \qquad E(e_3^2) = f(\lambda_{04} - 1),$$

$$E(e_1 e_p) = f\rho_{pb} C_\phi C_x, \qquad E(e_3 e_\phi) = fC_\phi \lambda_{12}, \qquad E(e_1 e_3) = fC_x \lambda_{03},$$

Where,

$$f = \left(\frac{1}{n} - \frac{1}{N}\right), \qquad C_p^2 = \frac{S_p^2}{P^2}, \qquad C_x^2 = \frac{S_x^2}{\bar{X}^2},$$

And $\rho_{pb}$ is the point bi-serial correlation coefficient.

## 2. Estimators in Literature

Singh et al. (2010) proposed the following ratio-type estimator for estimating unknown population mean when study variable is an attribute, as

$$t_a = \left(\frac{P}{\bar{x}}\right)\bar{X} \tag{2.1}$$

The bias and MSE expressions of the estimator $t_1$, to the first order of approximation is respectively, given by

$$B(t_a) = f\left(\frac{C_x^2}{2} - \rho_{pb}C_pC_x\right) \tag{2.2}$$

$$MSE(t_a) = fP^2\left(C_p^2 + C_x^2 - 2\rho_{pb}C_pC_x\right) \tag{2.3}$$

Singh et al. (2010) proposed a general class of estimator as,

$$t_b = H(p, u) \tag{2.4}$$

where $u = \dfrac{\bar{x}}{\bar{X}}$ and $H(p, u)$ is a parametric equation of p and u such that

$$H(p,1) = P, \forall P \tag{2.5}$$

and satisfying following regulations:

(i) Whatever be the sample chosen, the point (p,u) assume values in a bounded closed convex subset $R_2$ of the two-dimensional real space containing the point (p,1).

(ii) The function H(p,u) is a continuous and bounded in $R_2$.

(iii) The first and second order partial derivatives of H(p,u) exist and are continuous as well as bounded in $R_2$.

Where,

$$H_1 = \left.\frac{\partial H}{\partial u}\right|_{p=P, u=1}, \qquad H_2 = \left.\frac{1}{2}\frac{\partial^2 H}{\partial u^2}\right|_{p=P, u=1},$$

$$H_3 = \left.\frac{1}{2}\frac{\partial^2 H}{\partial p \partial u}\right|_{p=P, u=1}, \quad \text{and} \quad H_4 = \left.\frac{1}{2}\frac{\partial^2 H}{\partial p \bar{y}^2}\right|_{p=P, u=1}.$$

The bias and minimum MSE of the estimator $t_b$ are respectively, given by –

$$B(t_b) = f\left(P\rho_{pb}C_p C_x H_3 + C_x^2 H_2 + P^2 C_y^2 H_4\right) \tag{2.6}$$

$$MSE(t_b)_{min} = fP^2 C_p^2 \left(1 - \rho_{pb}^2\right) \tag{2.7}$$

Singh et al. (2010) proposed another family of estimator for estimating P, as

$$t_c = \left[q_1 P + q_2(\bar{X} - \bar{x})\right]\left[\frac{a\bar{X} + b}{a\bar{x} + b}\right]^\alpha \exp\left[\frac{(a\bar{X} + b) - (a\bar{x} + b)}{(a\bar{X} + b) + (a\bar{x} + b)}\right]^\beta \tag{2.8}$$

The bias and minimum MSE of the estimator $t_c$ to the first order of approximation, are respectively, given as

$$\text{Bias}(t_c) = P(q - 1) + f\left[(q_2 \bar{X} B + q_1 PA)C_x^2 - q_1 PB\rho C_p C_x\right] \tag{2.9}$$

$$MSE(t_c)_{min} = \left[P^2 - \frac{\Delta_1 \Delta_5^2 + \Delta_3 \Delta_4^2 - 2\Delta_2 \Delta_4 \Delta_5}{\Delta_1 \Delta_3 - \Delta_2^2}\right] \tag{2.10}$$

Where,

$$M_1 = P^2 f\left(C_p^2 + B^2 C_x^2 - 2B\rho C_p C_x\right), \qquad M_2 = \bar{X}^2 f\left(C_x^2\right),$$

$$M_3 = P^2 f\left(AC_x^2 - 2B\rho C_p C_x\right), \qquad M_4 = P\bar{X}f\left(-BC_x^2 + \rho C_p C_x\right),$$

$$M_5 = \overline{X}Pf(-BC_x^2)$$

$$q_1^* = \frac{\Delta_1\Delta_4 - \Delta_2\Delta_5}{\Delta_1\Delta_3 - \Delta_2^2} \quad \text{And} \quad q_2^* = \frac{\Delta_1\Delta_5 - \Delta_2\Delta_4}{\Delta_1\Delta_3 - \Delta_2^2} \tag{2.11}$$

Where,

$$\Delta_1 = (P^2 + M_1 + 2M_3), \Delta_2 = (-M_4 - M_5), \Delta_3 = (M_2), \Delta_4 = (P^2 + M_3), \Delta_5 = (-M_5),$$

## 3. Proposed estimators

The following estimator is proposed

$$t_1 = p\left(\frac{\overline{x}}{\overline{X}}\right)^\alpha \left(\frac{s_x^2}{S_x^2}\right)^\beta \tag{3.1}$$

Where $\alpha$ and $\beta$ are suitably chosen constants to be determined such that MSE of the estimator $t_1$ is minimum.

The bias and MSE of the estimator $t_1$ to the first order of approximation are respectively, given by

$$\text{Bias}(t_1) = \alpha\left(\frac{\alpha+1}{2}\right)C_p^2 + \beta\left(\frac{\beta+1}{2}\right)C_x^2 + \alpha\beta C_x\lambda_{03} + \alpha\rho_{px}C_pC_x + \beta C_p\lambda_{12} \tag{3.2}$$

$$\text{MSE}(t_1) = P^2f\left[C_p^2 + \alpha^2 C_x^2 + \beta^2(\lambda_{04}-1) + 2\alpha\rho_{px}C_pC_x + 2\beta C_p\lambda_{12} + 2\alpha\beta C_x\lambda_{03}\right] \tag{3.3}$$

Differentiating equation (3.3) partially with respect to $\alpha$ and $\beta$, equating them to zero, we get the optimum values of $\alpha$ and $\beta$ respectively, as

$$\alpha^* = \frac{C_p\{\lambda_{03}\lambda_{12} - \rho_{px}(\lambda_{04}-1)\}}{C_x\{(\lambda_{04}-1) - \lambda_{03}^2\}} \quad \beta^* = \frac{C_p\{\rho_{px}\lambda_{03} - \lambda_{12}\}}{\{(\lambda_{04}-1) - \lambda_{03}^2\}} \tag{3.4}$$

Putting the optimum values of α and β from equation (3.4) into equation (3.3), we get the minimum MSE of the estimator $t_1$ as

$$\text{MSE}(t_1)_{\min} = P^2 f C_p^2 \left[ 1 - \rho_{px}^2 - \frac{(\lambda_{03}\rho_{px} - \lambda_{12})^2}{(\lambda_{04} - 1 - \lambda_{03}^2)} \right] \tag{3.5}$$

Following Srivastava and Jhajj (1989), we propose a general family of estimators for estimating P, as

$$t_2 = H(p, u, v) \tag{3.6}$$

Where $u = \dfrac{\bar{x}}{\bar{X}}$, $v = \dfrac{s_x^2}{S_x^2}$ and $H(p, u, v)$ is a parametric function of p, u and v such that

$$H(p,1,1) = P, \forall P \tag{3.7}$$

And satisfying following regulations:

(iv) Whatever be the sample chosen, the point (p,u,v) assumes the values in a closed convex subset $R_3$ of the three-dimensional real space containing the point (p,1,1).

(v) The function H(p,u,v) is a continuous and bounded in $R_3$.

(vi) The first and second order partial derivatives of H(p,u,v) exist and are continuous as well as bounded in $R_3$.

Expanding H(p,u,v) about the point (P,1,1) in a second order Taylor series we have

$$t_2 = H(p, u, v) = [P + (p-P), 1 + (u-1), 1 + (v-1)] \tag{3.8}$$

$$(t_2 - P) = \left[ Pe_0 + e_1 H_1 + e_3 H_2 + Pe_0^2 H_3 + e_1^2 H_4 + e_3^2 H_5 + Pe_0 e_1 H_6 + e_1 e_3 H_7 + Pe_0 e_3 H_8 \right]$$

$$\tag{3.9}$$

Where,

$$\left.\frac{\partial H}{\partial p}\right|_{p=P, u=1} = 1,$$

$$H_1 = \left.\frac{\partial H}{\partial u}\right|_{p=P, u=1}, \qquad H_2 = \left.\frac{\partial H}{\partial v}\right|_{p=P, u=1, v=1},$$

$$H_3 = \left.\frac{1}{2}\frac{\partial^2 H}{\partial p^2}\right|_{p=P, u=1, v=1}, \qquad H_4 = \left.\frac{1}{2}\frac{\partial^2 H}{\partial u^2}\right|_{p=P, u=1, v=1}.$$

$$H_5 = \left.\frac{1}{2}\frac{\partial^2 H}{\partial v^2}\right|_{p=P, u=1, v=1}, \qquad H_6 = \left.\frac{1}{2}\frac{\partial^2 H}{\partial p \partial u}\right|_{p=P, u=1, v=1}$$

$$H_7 = \left.\frac{1}{2}\frac{\partial^2 H}{\partial u \partial v}\right|_{p=P, u=1, v=1}, \qquad H_8 = \left.\frac{1}{2}\frac{\partial^2 H}{\partial p \partial v}\right|_{p=P, u=1, v=1}$$

Taking expectations of both sides of (3.9), we get the bias of the estimator $t_2$ to the first order of approximation, as

$$B(t_2) = f\left[PC_p^2 H_3 + C_x^2 H_4 + (\lambda_{04} - 1)H_5 + P\rho_{px}C_p C_x H_6 + C_x \lambda_{03} H_7 + PC_p \lambda_{12} H_8\right] \quad (3.10)$$

Squaring both sides of (3.9) and neglecting terms of e's having power greater than two, we have

$$(t_2 - P)^2 = \left[Pe_0^2 + e_1^2 H_1^2 + e_3^2 H_2^2 + 2Pe_0 e_1 H_1 + 2Pe_0 e_3 H_2 + 2e_1 e_3 H_1 H_2\right] \quad (3.11)$$

Taking expectations of both sides of (3.11), we get the MSE of the wider class of estimator $t_2$ as

$$MSE(t_2) = f\left[PC_p^2 + C_x^2 H_1^2 + (\lambda_{04} - 1)H_2^2 + 2P\rho_{px}C_p C_x H_1 + 2PC_p \lambda_{12} H_2 + 2C_x \lambda_{03} H_1 H_2\right]$$

(3.12)

On differentiating (3.12) with respect to $H_1$ and $H_2$ equating to zero, respectively we obtain the optimum values of $H_1$ and $H_2$, as

$$H_1^* = \frac{C_p\{\lambda_{03}\lambda_{12} - \rho_{px}(\lambda_{04}-1)\}}{C_x\{(\lambda_{04}-1)-\lambda_{03}^2\}} \qquad H_2^* = \frac{C_p\{\rho_{px}\lambda_{03} - \lambda_{12}\}}{\{(\lambda_{04}-1)-\lambda_{03}^2\}} \qquad (3.13)$$

On substituting the values of $H_1^*$ and $H_2^*$ from (3.13) in (3.12), we obtain the minimum MSE of the estimator $t_2$, as

$$\text{MSE}(t_2)_{\min} = P^2 f C_p^2 \left[1 - \rho_{px}^2 - \frac{(\lambda_{03}\rho_{px} - \lambda_{12})^2}{(\lambda_{04} - 1 - \lambda_{03}^2)}\right] \qquad (3.14)$$

We propose another improved family of estimators for estimating P, as

$$t_3 = m_1 p \left[\frac{\overline{X}}{\gamma\overline{x} + (1-\gamma)\overline{X}}\right]^g + m_2 p \exp\left[\frac{\delta(S_x^2 - s_x^2)}{(S_x^2 + s_x^2)}\right] \qquad (3.15)$$

where $\gamma$ is suitable constant. g and $\delta$ are constants that can takes values (0,1,-1) for designing different estimators; and $m_1$ and $m_2$ are suitable chosen constants to be determined such that mean square error (MSE) of the class of estimator $t_3$ is minimum.

Expressing the class of estimators $t_3$ at equation (3.15) in terms of $e$'s, we have

$$t_3 = m_1 P(1+e_0)(1+\gamma e_1)^{-g} + m_2 P(1+e_0)\exp\left[\frac{-\delta e_3}{2}\left(1+\frac{e_3}{2}\right)^{-1}\right] \qquad (3.16)$$

Simplifying equation (3.16) and retaining terms up to the first order of approximation, we have

$$t_3 - P = -P\left[1 - m_1\left(1 - \gamma g e_0 e_1 + \frac{g(g+1)}{2}\gamma^2 e_1^2\right) - m_2\left(1 - \frac{\gamma}{2}e_0 e_1 + \frac{\delta(\delta+1)}{8}e_1^2\right)\right]$$

(3.17)

Taking expectations of both sides of equation (3.17), we get the bias of the estimator $t_3$ to the first order of approximation, as

$$\text{Bias}(t_3) = -P[1 - m_1 B - m_2 E] \tag{3.18}$$

where, $B = \left\{1 - \gamma g f \rho_{px} C_p C_x + \dfrac{g(g+1)}{2}\gamma^2 f C_x^2\right\}$

$$E = \left\{1 - \dfrac{\gamma}{2} f \rho_{px} C_p C_x + \dfrac{\delta(\delta+1)}{8} f C_x^2\right\} \tag{3.19}$$

Squaring both sides of equation (3.17) and neglecting the terms having power greater than two, and then taking expectations of both sides, we get the MSE of the estimator $t_3$ to the first order of approximation, as

$$\text{MSE}(t_3) = P^2[1 + m_1^2 A + m_2^2 C + 2m_1 m_2 D - 2m_1 B - 2m_2 E] \tag{3.20}$$

Where, $A = \left\{1 + f\left(C_p^2 - 4\gamma g \rho_{px} C_p C_x + \gamma^2 g(2g+1) C_x^2\right)\right\}$

$$C = \left\{1 + f\left(C_p^2 - 2\delta C_p \lambda_{12} + \dfrac{\{\delta^2 + \delta(\delta+2)\}(\lambda_{04} - 1)}{4}\right)\right\}$$

$$D = \left\{1 + f\left(C_p^2 - \delta C_p \lambda_{12} + \dfrac{\delta(\delta+2)}{8}(\lambda_{04} - 1) - 2\gamma g \rho_{px} C_p C_x \right.\right.$$

$$\left.\left. + \dfrac{\gamma \delta g}{2} C_x \lambda_{03} + \dfrac{g(g+1)}{2}\gamma^2 C_x^2 \right)\right\}$$

And B and E are the same as defined earlier give in (3.19).

The MSE of the class of estimator $t_3$ at equation (3.20) is minimised for the optimum values of $m_1$ and $m_2$ given as

$$m_1^* = \dfrac{(BC - DE)}{(AC - D^2)} \quad \text{and} \quad m_2^* = \dfrac{(AE - BD)}{(AC - D^2)} \tag{3.21}$$

Putting equations (3.21) in (3.20), we get the resulting minimum bias and MSE of the proposed class of estimators $t_3$, respectively, as

$$\text{Bias}(t_3)_{\min} = -P\left[1 - \dfrac{B^2 C - 2BDE + AE^2}{(AC - D^2)}\right] \tag{3.22}$$

$$\text{MSE}(t_3) = P^2 \left[ 1 - \frac{B^2 C - 2BDE + AE^2}{(AC - D^2)} \right] \tag{3.23}$$

## 4. Efficiency Comparisons

First we compare the MSE of proposed estimators $t_1$ and $t_2$ with usual estimator,

$$\text{MSE}(t_1)_{\min} = \text{MSE}(t_2)_{\min} \leq V(p)$$

If,

$$fC_p^2 \left[ 1 - \rho_{px}^2 - \frac{(\lambda_{03}\rho_{px} - \lambda_{12})^2}{(\lambda_{04} - 1 - \lambda_{03}^2)} \right] \leq fC_p^2 \tag{4.1}$$

On solving we observed that above conditions holds always true.

Now, we compare the efficiency of proposed estimator $t_3$ with usual estimator,

$$\text{MSE}(t_3)_{\min} \leq V(p)$$

If,

$$\left[ 1 - \frac{B^2 C - 2BDE + AE^2}{(AC - D^2)} \right] \leq f_1 C_p^2 \tag{4.2}$$

On solving we observed that above conditions holds always true.

Next we compare the efficiency of proposed estimator $t_3$ with wider class of estimator $t_2$.

$$\text{MSE}(t_3)_{\min} \leq \text{MSE}(t_2)_{\min} = \text{MSE}(t_1)_{\min}$$

If,

$$\left[ 1 - \frac{B^2 C - 2BDE + AE^2}{(AC - D^2)} \right] \leq fC_p^2 \left[ 1 - \rho_{px}^2 - \frac{(\lambda_{03}\rho_{px} - \lambda_{12})^2}{(\lambda_{04} - 1 - \lambda_{03}^2)} \right] \tag{4.3}$$

Finally, we compare the efficiency of proposed estimator $t_3$ with class of estimator $t_c$ proposed by Singh et al. (2010) as

$$MSE(t_3)_{min} \leq MSE(t_c)_{min}$$

if

$$P^2\left[1 - \frac{B^2C - 2BDE + AE^2}{(AC - D^2)}\right] \leq \left[P^2 - \frac{\Delta_1\Delta_5^2 + \Delta_3\Delta_4^2 - 2\Delta_2\Delta_4\Delta_5}{\Delta_1\Delta_3 - \Delta_2^2}\right] \quad (4.4)$$

## 5. Empirical study

**Data Statistics:** The data used for empirical study has been taken from Gujrati and Sangeetha (2007) -pg, 601.

Here,

Y – Home ownership.

X – Income (thousands of dollars).

| n | N | P | $\overline{X}$ | $\rho_{pb}$ | Cp | Cx | $\lambda_{12}$ | $\lambda_{04}$ | $\lambda_{03}$ |
|---|---|---|---|---|---|---|---|---|---|
| 11 | 40 | 0.525 | 14.4 | 0.897 | 0.963 | 0.308 | -0.118 | 1.75 | -0.153 |

The following Table shows comparison between some existing estimators and proposed estimators with respect to usual estimator.

**Table 5.1: Percent relative efficiency of proposed estimators with respect to usual estimator**

| Estimators | p | $t_a$ | $t_b$ | $t_c$ | $t_1$ | $t_2$ | $t_3$ | | |
|---|---|---|---|---|---|---|---|---|---|
| | | | | | | | $g=1, \delta=1$ | $g=1, \delta=-1$ | $g=0, \delta=1$ |
| PRE | 100 | 189.38 | 511.79 | 518.05 | 513.92 | 513.92 | **685.51** | 199.20 | 141.23 |

When we examine Table 5.1, we observe that the proposed estimators $t_1$, $t_2$ and $t_3$ all performs better than the usual estimator p. Also, the proposed estimator $t_3$ is the best among the estimators considered in this paper and perform better than the estimators proposed by Singh et al. (2010) for estimating P for the choice $g = 1, \delta = 1$.